\documentclass[12pt]{amsart}
\usepackage{amssymb}
\usepackage[dvips]{graphics}
\ifx\mathrm\undefined\let\mathrm\rm\fi
\ifx\mathbf\undefined\let\mathbf\bf\fi
\ifx\mathfrak\undefined\let\mathfrak\frak\fi
\ifx\mathcal\undefined\let\mathcal\cal\fi
\ifx\mathbb\undefined\let\mathbb\Bbb\fi
\ifx\emph\undefined\let\emph\it\fi
\font\bb=msbm10 at9.98pt
\begin{document}
\def\semidirect{\hbox{$\;$\bb\char'156$\;$}}
\newcommand{\SL}{\mathrm{SL}}
\newcommand{\GL}{\mathrm{GL}}
\newcommand{\g}{{{\mathfrak g}\,}}
\newcommand{\bor}{{{\mathfrak b}}}
\newcommand{\n}{{{\mathfrak n}}}
\newcommand{\h}{{{\mathfrak h\,}}}
\newcommand{\Id}{{\operatorname{Id}}}
\newcommand{\Z}{{\mathbb Z}}
\newcommand{\ZZ}{{\mathbb Z_{>0}}}
\newcommand{\N}{{\mathbb N}}
\newcommand{\R}{{\mathbb R}}
\newcommand{\p}{{\mathbb P}} 
\newcommand{\C}{{\mathbb C}}
\newcommand{\Q}{{\mathbb Q}}
\newcommand{\CC}{\mathcal{C}}
\newcommand{\A}{\mathcal{A}}
\newcommand{\F}{\mathcal{F}}
\newcommand{\W}{\mathcal{W}}      
\newcommand{\PP}{\mathcal{P}}
\newcommand{\Poly}{\rm{Poly}}
\newcommand{\Span}{\rm{Span}}
\newcommand{\Sing}{\rm{Sing}}
\newcommand{\Sym}{\rm{Sym}}
\newcommand{\1}{{\bf 1}}
\newcommand{\m}{{\bf m}}
\newcommand{\z}{{\bf z}}
\newcommand{\s}{{\bf s}}    
\newcommand{\tb}{{\bf t}}
\newcommand{\jj}{{\bf j}}                            
\newcommand{\dontprint}[1]
{\relax}
\newtheorem%
{thm}{Theorem}
\newtheorem%
{prop}
{Proposition}
\newtheorem%
{lemma}
{Lemma}
\newtheorem%
{lemmadef}[thm]{Lemma-Definition}
\newtheorem%
{cor}
{Corollary}
\newtheorem%
{conj}
{Conjecture}
\newenvironment{definition}
{\noindent{\bf Definition\/}:}{\par\medskip}
                                                    
\title {Rational functions with prescribed critical points}
 
\author[{}]
{I. Scherbak}

\maketitle

\medskip
\centerline{\it School of Mathematical Sciences,
Tel Aviv University,}
\centerline{\it Ramat Aviv, Tel Aviv 69978, Israel}                           
\centerline{\it e-mail: \quad scherbak@post.tau.ac.il}

\bigskip 

\pagestyle{myheadings}
\markboth{I. Scherbak}
{ Rational functions  with prescribed critical points}
\begin{abstract}
${}$
A rational function is the ratio of two complex polynomials in one
variable without common roots.
Its  degree is the maximum of the degrees of the numerator and
the denominator.  Rational functions belong to the same class
if one turns into the other by the postcomposition with a 
linear-fractional transformation. We give an explicit formula
for the number of classes having a given degree $d$ and given 
multiplicities $m_1,\dots,m_n$ of given $n$ critical points,
for generic positions of  the critical points.
This number is the multiplicity of the irreducible
$sl_2$ representation with highest weight
$2d-2-m_1-\dots -m_n\,$
in the tensor product of the irreducible
$sl_2$ representations with highest weights  $m_1,\dots,m_n$.

The classes are labeled by the orbits of critical points of a remarkable 
symmetric function which  first appeared in the XIX century in studies 
of Fuchsian differential equations, and then in  the XX century in 
the theory of KZ equations.
\end{abstract}

\section {Introduction}

{\it A rational function} is the ratio of two complex polynomials 
in one variable without common roots.
Its {\it degree} is the maximum of the degrees of the numerator and 
the denominator. Any rational function defines a branched covering 
of the Riemann sphere, $\C\p^1=\C\cup\infty$, and the degree is 
the cardinality of the preimage of any regular value. 

\medskip
Let $R(x)=g(x)/f(x)$, where $g(x)$ and $f(x)$ are polynomials
without common roots, be a rational function of degree $d$.
We have $R'(x)=W(x)/f^2(x)$, where 
$$
W(x)=W[g,f](x)=g'(x)f(x)-g(x)f'(x)
$$
is the {\it Wronskian of polynomials $g$ and $f$}.

\medskip
The rational function  $ R(x)$  
has {\it a critical point of multiplicity $m\geq 1$}  
at  $z\in\C$, if $z$ is a root of the Wronskian
$W(x)$ of multiplicity $m$.
Let $m_1,\dots,m_n$ be positive integers and $z_1,\dots, z_n$			
pairwise distinct complex numbers.  Write $\m=(m_1,\dots,m_n)$,			
$M=m_1+\dots+m_n\,$ and $\z=(z_1,\dots, z_n)$.					
If $R(x)$ has degree $d>1$ and if all its finite 				
critical points are $z_1,\dots, z_n$  of multiplicities  			
$m_1,\dots,m_n$ respectively, then
we say that $R(x)$ {\it has the type $(d,n;\m;\z)$}.
According to the Riemann-Hurwitz formula, numbers 
$d$ and $M$ should satisfy $2d-2\geq  M$. The difference 
$m_{\infty}=2d-2-M$ is {\it the multiplicity of  $R(x)$ 
at infinity}.
					
\medskip
If $R(x)$ has the type $(d,n;\m;\z)$, 
then the rational function given by the postcomposition with 
a linear-fractional transformation,
$$ 
\frac{aR(x)+b}{cR(x)+d}\,,\quad ad-bc\neq 0\,,
$$ 
has clearly the same type. 
A rational function considered up to postcompositions
with linear-fractional transformations will be called 
{\it a class of rational functions}. 

\medskip
We address the following question.

\medskip\centerline
{\it  Given $d,n,\m,\z$, how many classes have the type $(d,n;\m;\z)$?}

\medskip
The question  is easy  for   $n=1$.
If  $z\in\C$  is the single critical point,  then the number of 
classes is $1$ if  the multiplicity of $z$ is $d-1$ (this is the class containing  
$(x-z)^d$), and $0$ otherwise.

If  $n\geq 2$, we answer the question for generic $\z$. 
The words ``$\z$ is  generic'' mean that $\z$ does not belong to a suitable proper 
algebraic set in $\C^n$. 

\bigskip\noindent{\bf Theorem.}\quad
{\it Let  $n\geq 2$ and $ d,\, m_1\,, \dots, m_n\,$  be positive integers.
For generic $\z$,  the number $\sharp(d,n;\m)$ of  classes of rational functions 
of the type $(d,n;\m;\z)$ is}
\begin{eqnarray}\label{sharp}
\sharp(d,n;\m)=\sum_{q=1}^n (-1)^{n-q}
\sum_{1\leq i_1<\dots <i_q\leq n}
{m_{i_1}+\dots +m_{i_q}+q-d-1\choose n-2}\,,
\end{eqnarray}
{\it and  any nonempty class can be represented by the ratio 
of polynomials without multiple roots.}

\medskip\noindent
As usually we set ${a\choose b}=0$ for $a<b$. 

\medskip\noindent
{\bf Remarks:}\ \
\begin{itemize}
\item
If   $m_j>d-1$ for some $j$, or if  $M>2d-2$, or if
$M<d-1$,  then the right hand side of (\ref{sharp}) vanishes.
In fact, in this case for  {\it any} $\z$ there are no 
rational functions of the type $(d,n;\m;\z)$. Indeed, 
any of the inequalities is impossible for  a rational map of degree $d$. 
This was known,  see Ch.~6 of \cite{GH} or  Lemma~1.1 of \cite{G}.
 \item
If  $M=d-1$ and   $m_j\leq d-1$ for  $j=1,\dots,n$, then  
the right hand side of (\ref{sharp}) gives $1$.  Again, for {\it any} $\z$ 
the number of classes is $1$.  Indeed, if  $R(x)=g(x)/f(x)$ belongs to a 
corresponding  class,  and if $g$ has degree $d$, then the  degree of $f$ is 
clearly zero,  i.e. this class contains the polynomial 
$$
\int (x-z_1)^{m_1}\dots (x-z_n)^{m_n}\, dx\,.
$$
\item  
The case when the number of critical points is the maximal possible, i.e.
 $n=2d-2$ and $m_1=\dots=m_n=1$, was studied  by L.~Goldberg
 in  \cite{G}. Denote  $\1=(1,\dots,1)$.  Theorem~1.3  of    \cite{G} says
that  $\sharp(d,2d-2;\1)=\frac 1d{2d-2\choose d-1}\,$,
which is the $d$-th Catalan number.
\item
The condition ``$\z$ is generic''  is essential.
For  some of $\z$, the number of classes may  decrease. 
Actually, the inequality (\ref{s1}) and the left-hand inequality of (\ref{s3})
hold for any $\z$ whereas the right-hand inequality of (\ref{s3}) may be broken
for  some of $\z$. An example is given in Theorem 1.4 of  \cite{G}, cf. Sec.~12 
of \cite{RV}.
\end{itemize}

In view of the first and second remarks, we assume
\begin{equation}\label{A}
d-1 < M\leq 2d-2,\quad 1\leq m_j\leq d-1,\quad j=1,\dots,n.
\end{equation}
The proof of the theorem consists of the following three steps.

\bigskip\noindent \underline{{\it STEP  1.}}\quad  
By means of the Schubert calculus
we obtain an upper bound  by   the intersection number  
of  special  Schubert 
classes $\sigma_{m_1}$, ..., $\sigma_{m_n}$, $\sigma_{2d-2-M}$:

\begin{equation}\label{s1}
\sharp(d,n;\m) \leq
\sigma_{m_1}\,\cdot {\rm\ ...\ }\cdot \sigma_{m_n}\,\cdot\, \sigma_{2d-2-M}\,.
\end{equation}
This step is done in Sec.~\ref{S2}.

\bigskip\noindent\underline{ {\it STEP  2.}}\quad  Denote  $L_{m_j}$  the irreducible 
$sl_2$-module with  highest weight $m_j$, and   ${\Sing}_k L$  the subspace 
of  singular vectors of weight  $M-2k$ in the tensor product 
$L=L_{m_1}\otimes \dots \otimes  L_{m_n}$. Well-known  relations 
of the Schubert  calculus to the representation theory (\cite{F}) imply
\begin{equation}\label{s2}
\sigma_{m_1}\,\cdot {\rm\ ...\ }\cdot \sigma_{m_n}\,\cdot\, \sigma_{2d-2-M}\,
=\, \dim{\Sing}_{M+1-d} L.
\end{equation}
The number  $ \dim{\Sing}_k L$ was calculated in Theorem~5 of \cite{SV}:
\begin{equation}\label{dim}
\dim{\Sing}_kL\,=\,\sum_{q=0}^{n} (-1)^q
\sum_{1\leq i_1<\dots <i_q\leq n}
{k+n-2-m_{i_1}-\dots -m_{i_q}-q\choose n-2}.
\end{equation}

\bigskip
In order to replace the inequality  (\ref{s1}) with 
the equality  one needs a non-trivial transversality statement which
we did not find  in the algebraic geometry literature 
\footnote
{As  the author learned  from the referee report, 
after the first version of this paper 
was posted on  arXiv and submitted to this journal, a preprint of  B.~Osserman
has appeared.  The preprint contains an algebraic geometry proof  
of the equality 
$\sharp(d,n;\m) =\sigma_{m_1}\,\cdot {\rm\ ...\ }\cdot \sigma_{m_n}\,\cdot\, 
\sigma_{2d-2-M}\,$ based on the results  
of  D.~Eisenbud and J.~Harris, \cite{EH}.}.

We extract  a lower bound  from  the results of  
N.~Reshetikhin and A.~Varchenko   
on  the $sl_2$  Knizhnik--Zamolodchikov  (KZ)  equation, \cite{RV}. 
The key observation is that rational functions and   
asymptotic solutions to  the  
$sl_2$ KZ equations  turn to be linked by a remarkable  symmetric {\it master
function} whose critical  point orbits  determine  both classes of rational 
functions and asymptotic solutions to  the  $sl_2$ KZ equations. 
This  is the very  function  studied by  Heine and Stieltjes 
in connection with second order Fuchsian  differential equations having 
a polynomial solution of a prescribed degree.   

\bigskip\noindent\underline{{\it STEP  3.}}\quad 
Denote \{NCO\} the number of orbits of critical points of the master 
function given by (\ref{Phi}). We get
\begin{equation}\label{s3}
\sharp(d,n;\m) \geq \{{\rm NCO}\}\geq  \dim{\Sing}_{M+1-d} L.
\end{equation}
The left-hand inequality follows from  the classical result of
Heine and Stieltjes, Ch.~6.8 of  \cite{Sz}. 
We describe the connection between rational functions and 
Fuchsian equations and prove this inequality in  Sec.~\ref{S3}.
The right-hand inequality is proved in Theorem~9.9 of  \cite{RV} 
(cf. Theorem~8 of  \cite{SV}). The appearance of the master function 
in the theory of  KZ equations is explained in Sec.~\ref{S4}, 
following  \cite{RV}.

\medskip
Juxtaposing relations (\ref{s1})--(\ref{s3}) we arrive at the result.

\medskip
Sec.~\ref{S5} contains comments. In particular, connections of the
Schubert calculus with Fuchsian equations having only polynomial solutions
(Sec.~\ref{s54}) and with  $sl_p$ KZ equations (Sec.~\ref{s55}) are discussed.  

 \medskip\noindent
{\bf Acknowledgments.}\ \
The author is grateful to  A.~Eremenko, A.~Gabrielov, B.~Shapiro, 
M.~Sodin, F.~Sottile, A.~Vainshtein, A.~Varchenko,  V.~Zakalyukin,
and the referee  for useful comments and discussions.

\section{An upper bound:  Schubert classes}\label{S2}
\subsection{Classes of rational functions and generic planes in the
Grassmannian}\label{s21}

Denote  $G_2(\Poly)$  the Grassmannian of  two-dimensional 
 planes in  the vector space of  complex polynomials
 in one variable. Any plane $V\in G_2(\Poly)$ has a basis of polynomials
of distinct degrees.  Define  {\it the degree}  of $V$ as  the maximal 
degree of   polynomials in $V$, and  {\it the order} of $V$ 
as  the minimal  degree of  non-zero polynomials in $V$.
Say that $V\in G_2(\Poly)$  is {\it a generic plane}, if
for any  $z\in\C$ there exists a polynomial $p(x)\in V$ such that 
$p(z)\neq 0$. 

\begin{prop}\label{p1} There is a one-to-one correspondence
between the generic planes and the classes of rational
functions of the same degree.
\end{prop}

\noindent
{\bf Proof:}\quad   $V\in G_2(\Poly)$ is generic if and only if
 the  polynomials of any basis of $V$ do not have common roots.
 Any rational function  $R(x)=g(x)/f(x)$ clearly defines a
generic plane $V_R= {\rm Span}\{g(x), f(x)\}$ of the same degree.
Conversely, any basis $\{g(x), f(x)\}$ of 
a generic plane $V$ of degree $d$ defines a rational function  
$R_V(x)=g(x)/f(x)$  of  degree $d$.
Linear-fractional transformations of  the rational function 
correspond to changes of  the basis.
\hfill $\triangleleft$

\subsection{The Wronski map}\label{s22}
 Define {\it the Wronskian of} $V\in  G_2(\Poly)$  as a monic
 polynomial  which is proportional to the Wronskian
 of  some (and hence, any) basis of $V$.  The question
 posed in the Introduction can be reformulated as follows:

 \medskip\centerline
 {\it Given monic polynomial $W(x)$ and  given positive integer $d$,}

\centerline
 {\it how many generic planes of degree $d$  have the
Wronskian $W(x)$?}

\medskip
 It is enough to consider the  $(d+1)$-dimensional vector space ${\Poly}_d$
of polynomials of degree $\leq d$ and the Grassmannian 
$G_2({\Poly}_d)$ of two-dimensional subspaces in ${\Poly}_d$.
We have $\dim G_2({\Poly}_d)=2d-2$.
The maximal possible degree of the Wronskian of an element of  
$G_2({\Poly}_d)$  is $2d-2$.  Polynomials of degree $\leq 2d-2$, considered up
to a non-zero factor, form a $(2d-2)$-dimensional complex projective space.
 Thus the mapping  of $V\in G_2({\Poly}_d)$ into the Wronskian is  a well-defined 
 map  of   smooth complex algebraic varieties  of the same dimension,
 $G_2({\Poly}_d)\to \C\p^{2d-2}$, called  {\it the Wronski map}. The preimage of  
 any monic polynomial consists of  a finite number of   planes.  
 On Wronski maps see \cite{EG}. 

\subsection{Schubert calculus {\rm (Ch.~6 of \cite{GH})}}\label{s23}
Let   $G_2(d+1)$ be the Grassmannian variety of  
two-dimensional subspaces $V\subset \C^{d+1}$.   
A chosen basis  $e_1,\dots,e_{d+1}\,$ of $ \C^{d+1}$
defines the flag of linear subspaces
$$
E_{\bullet}\,:\quad E_1\, \subset \, E_2\, \subset\, \dots\,  
\subset \, E_d \subset \, E_{d+1}= \C^{d+1}\,,
$$
where $E_i={\Span}\{e_1,\dots,e_i\}\,$, $\dim E_i=i$.   
For any integers  $a_1$ and $a_2$
such that   $0\leq a_2\leq a_1\leq d-1\,$, the {\it Schubert variety} 
$\Omega_{a_1,a_2}(E_{\bullet})\subset  G_2({d+1})$
is defined as follows,
$$
\Omega_{a_1,a_2}= \Omega_{a_1,a_2}(E_{\bullet})=\left\{\, V\in  G_2({d+1})\, 
\vert\, \dim \left(V\cap E_{d-a_1}\right)\geq 1\,,
\   \dim \left(V\cap E_{d+1-a_2}\right)\geq 2\, \right\}\,.
$$
The variety $\Omega_{a_1,a_2}=\Omega_{a_1,a_2}(E_{\bullet})$
is an irreducible closed subvariety of $ G_2({d+1})$  of the complex
codimension $a_1+a_2$.

\medskip
The homology classes $[\Omega_{a_1,a_2}]$  of Schubert varieties
 $\Omega_{a_1,a_2}$ are independent of the choice of  flag, and form a basis
for the integral homology of  $G_2(d+1)$.
Define  $\sigma_{a_1,a_2}$ to be the cohomology class in 
 $H^{2(a_1+a_2)}(G_2({d+1}))$ whose cap product with the
fundamental class of   $G_2(d+1)$ is the homology class
 $[\Omega_{a_1,a_2}]$.
The classes  $\sigma_{a_1,a_2}$ are called {\it Schubert classes}.
They give a basis over $\Z$ for the cohomology ring of the Grassmannian.
The product or {\it intersection} of  any two Schubert classes   
$\sigma_{a_1, a_2}$ and  $\sigma_{b_1, b_2}$
has the form
$$
\sigma_{a_1, a_2}\,\cdot\,\sigma_{b_1, b_2}\,=\,
\sum_{c_1+c_2=a_1+a_2+b_1+b_2} C(a_1,a_2;b_1,b_2; c_1,c_2) 
\sigma_{c_1, c_2}\,,
$$
where $ C(a_1,a_2;b_1,b_2; c_1,c_2)$ are nonnegative integers
called {\it the Littlewood--Richardson coefficients}.

\medskip
If the sum of the codimensions of classes equals $\dim G_2({d+1})=2d-2$,
then their intersection is an integer (identifying the generator of the top
cohomology group $\sigma_{d-1, d-1}\in H^{4d-4}(G_2({d+1}))$ with $1\in\Z$) 
called {\it the intersection number}.

\medskip
When $(a_1, a_2)=(q,0)$, $0\leq q\leq d-1$, 
the  Schubert varieties $\Omega_{q, 0}$  
are called {\it special} and the corresponding cohomology classes 
$\sigma_q\,=\,\sigma_{q, 0}$
are called {\it special Schubert classes}.

 \subsection{Planes with a given Wronskian}\label{s24}
 In order to apply the Schubert calculus to our problem, we need
 the following result.

\begin{prop}\label{p2}  Let $V\in G_2({\Poly}_d)$  be  a generic plane 
of degree $d$ and order $k$ with the Wronskian $W(x)$.  Then

\begin{itemize}
\item
almost all polynomials  of  degree $d$ in $V$ do not have multiple
roots;
\item
the  polynomials of  degree $k$ in $V$ are all proportional;
 \item
if there exists a polynomial  $p(x)\in V$  which has a root 
of multiplicity $1$ at $0$,  then  $W(0)\neq 0$ and 
$$
\frac{W'(0)}{W(0)}\,=\,\frac{p''(0)}{p'(0)}\,.
$$  
\item
if   $0$  is  a root of multiplicity  $m\geq 1$ of   $W(0)$,
  then there exists a polynomial  in $V$ which has a root
  of multiplicity  $m+1$ at $0$.  
\end{itemize}
\end{prop}

\noindent
{\bf Proof:}\ \   The first two statements   are evident. 

In order to prove the third statement, take  a linearly  independent 
polynomial $q(x)\in V$.
We  have $q(0)\neq 0$ and $W[p,q](x)=p'(x)q(x)-p(x)q'(x)=cW(x)$ for 
some constant $c\neq 0$. The substitution of $x=0$ gives 
$cW(0)=p'(0)q(0)\neq 0$. 
Moreover $cW'(x)=p''(x)q(x)-p(x)q''(x)$, and the substitution of $x=0$ into
$$
\frac{W'(x)}{W(x)}\,=\,\frac{p''(x)q(x)-p(x)q''(x)}{p'(x)q(x)-p(x)q'(x)}
$$ 
finishes the proof of the third statement.

In order to prove the last statement, take a basis of $V$,
$\{f(x)\,,\, g(x)\}$  such that  $\deg f <\deg g$.
We have  $W(x) =x^m\tilde W(x)$, where $\tilde W(0)\neq 0$.
There are two possibilities:  $f(0)=0$  or $f(0)\neq 0$.

\medskip
1) If $f(0)=0$, then  $g(0)\neq 0$ and  $f(x)$ has the form 
$f(x)=x^k\tilde f(x)$, 
where $\tilde f(0)\neq 0$,  for some $k >1$. The Wronskian of  $f$ and $g$,
which is proportional to $W(x)$, has the form $x^{k-1}F(x)$, where
$F(0)= k\tilde f(0)g(0)\neq 0$. Therefore $k=m+1$.

\medskip
2) If $f(0)\neq 0$, then 
$$
g_0(x)\,=\,g(x)-\frac{g(0)}{f(0)}\,f(x)\,\in V
$$
satisfies  $g_0(0)=0$, and   polynomials $f(x)\,,\, g_0(x)$  form
a basis of $V$. The rest of the proof is the same as in  case 1).
\hfill $\triangleleft$

\subsection{The proof of  {\rm (\ref{s1})}}\label{s25}
Now we are able   to estimate from above the number of generic 
planes of degree $d$ with the Wronskian
\begin{equation}
\label{W}
W(x)=W(x;\z,\m)=(x-z_1)^{m_1}\dots (x-z_n)^{m_n}\,, \quad \deg W=M.\\
\end{equation}

\medskip
For  any fixed $z\in\C$, define the  flag  $\F_\bullet(z)$ in $ {\Poly}_{d}\,$, 
\begin{equation*}
 \F_\bullet(z)\,:\,\  \F_0(z) \subset \F_1(z)\subset \dots 
\F_{d}(z)= {\Poly}_{d}\,,
 \quad \dim \F_i(z)= i+1\,,
\end{equation*}
where  $\F_i(z)$  consists of the polynomials  of the form
$$
 a_i(x-z)^{d-i}+\dots +a_{0}(x-z)^{d}\,.
 $$
Define $\F_\bullet(\infty)$ as the flag with  $\F_i(\infty)={\Poly}_i\,$, 
$0\leq i\leq d$.

\medskip
Let $V$ be a generic plane of degree $d$ with the Wronskian given by
(\ref{W}) and (\ref{A}).

\medskip
According to Proposition \ref{p2}, for any $z_j$ there exists a basis of $V$
such that one of the polynomials does not vanish at $z_j$ and  the other
has a root at $z_j$ of multiplicity $m_j+1$.  This exactly means
that  $V\in \Omega_{m_j, 0}(\F_{z_j})$. 

\medskip
The Wronskian of $V$ has degree $M$, therefore $V$ should 
contain a polynomial of degree $M+1-d$.  
Thus  $V$ has a basis of polynomials of degrees
$d$ and $M+1-d$, i.e.  $V\in  \Omega_{2d-2-M, 0}(\F_{\infty})$.

\medskip
Thus we conclude that
$$
V\in \Omega_{m_1, 0}(\F_\bullet(z_1))\,\cap \Omega_{m_2, 0}(\F_\bullet(z_2))
\cap \dots \cap \Omega_{m_n, 0}(\F_\bullet(z_n))
\, \cap \Omega_{2d-2-M, 0}(\F_\bullet({\infty}))\,.
$$
The dimension of  $G_2({\Poly}_{d})$  is the sum of the codimensions
of  these Schubert varieties, therefore  the intersection 
consists of a finite number of planes. This number does not exceed
$\sigma_{m_1}\,\cdot {\rm\ ...\ }\cdot \sigma_{m_n}\,\cdot\, 
\sigma_{2d-2-M}\,$.
\hfill $\triangleleft$     

\section{The master function}\label{S3}
\subsection{Fuchsian differential equations with only polynomial
solutions}\label{s31}
Consider a second order Fuchsian differential equation with 
regular singular points at $z_1,\dots,z_n$, $n\geq 2$, and at 
infinity.  If   the exponents at  $z_j\,$ are $0$ and $m_j+1$,
$1\leq j\leq n$, then this equation has  the form
\begin{eqnarray}\label{FE}
F(x)u''(x)+G(x)u'(x)+H(x)u(x)=0\,,\\
F(x) = \prod_{j=1}^n (x-z_j)\,,\quad
\frac{G(x)}{F(x)} = \sum_{j=1}^n\frac{-m_j}{x-z_j}\,,\nonumber
\end{eqnarray}
where $H(x)$ is a polynomial of degree not greater than $n-2$.
On Fuchsian equations see Ch.~6 of \cite{R}.
As before, we write  $M=m_1+\dots+m_n$.

\begin{prop}\label{pF}
Let $m_1,\dots,m_n$, $n\geq 2$, and $d$ be positive integers such that
 $d>M+1-d\geq 0$.
There is a one-to-one correspondence between the
generic planes of  degree $d$ with the Wronskian {\rm (\ref{W})}
and the equations  {\rm (\ref{FE})} with only polynomial solutions and the 
degree $d$ generic solutions.
\end{prop}
 
\noindent
{\bf Proof:}\quad 
The Wronskian  of  the equation  (\ref{FE})  is exactly the polynomial $W(x)$
given by (\ref{W}).  
If this equation has only polynomial solutions and if  the generic solution has 
degree $d$, then the solution space is a plane in  $G_2({\Poly}_d)$ of degree $d$.  
If  $g(x)$ and $f(x)$ form a basis in  the solution space, 
then  the Wronskian  of  $g$ and $f$  is proportional to  $W(x)$.  
Therefore  $g$  and $f$  may have a common root  at points $z_j$ only. 
But the generic solution to a Fuchsian equation   (\ref{FE}) can not vanish 
at a regular point (Ch.~6 of \cite{R}). 
 
 \medskip
In order to prove the opposite statement, consider a
basis  $\{g(x),f(x)\}$ of a generic plane $V$
of  degree $d$ with the Wronskian (\ref{W}). 
We can assume 
$f(x)=(x-t_1)\cdots(x-t_d)\,$ where the roots
$t_1,\dots, t_d$ satisfy
$$
t_i\neq t_l\,,\ t_i\neq z_j\,,\ 1\leq i\,, l\leq d\,,\
1\leq j\leq n\,,
$$
according to Proposition  \ref{p2}.
The plane $V$  is the solution space of the following second order   
linear differential equation with respect to unknown function $u(x)$,
$$
\left|\begin{array}{ccc}u(x)\ \ &f(x)\   \ &g(x)\\
                          u'(x)\ \ &f'(x)\ \  &g'(x)\\
                          u''(x)\  \ &f''(x)\ \ &g''(x)
                           \end{array}\right|=0.
$$
The Wronskian of the polynomials $f(x)$ and $g(x)$ is proportional
to $W(x)$ given by (\ref{W}), therefore this equation can be re-written 
in the form
$$
W(x)u''(x)-W'(x)u'(x)+h(x)u(x)=0\,,
$$
where $h(x)$ is a polynomial  proportional to the Wronskian of
$f'(x)$ and $g'(x)$. One can easily check that
$$
\frac{W'(x)}{W(x)}\,=\,\sum_{j=1}^n\frac{m_j}{x-z_j}\,.
$$
Moreover, we have
$$
h(x)\, =\, \frac{-W(x)f''(x)+W'(x)f'(x)}{f(x)}\,,
$$
as $f(x)$ is clearly  a solution to this equation.
Hence if $z_j$ is a root of $W(x)$ of multiplicity $m_j>1$,
then all coefficients of the equation have $(x-z_j)^{m_j-1}$ as
a common factor, and the equation can be reduced to the required 
form (\ref{FE}).  
\hfill $\triangleleft$

\begin{cor}\label{c2}
Any  plane in $G_2(\Poly)$ with a given Wronskian is uniquely determined by any 
of its  polynomial.            
\hfill $\triangleleft$    
\end{cor}

\medskip
The following statement  easily follows from the theory of Fuchsian equations 
(Ch.~6 of \cite{R}).

\begin{prop}\label{pp} {\rm (Cf. Sec.~3.1 of \cite{SV})}\ \
If the equation  {\rm  (\ref{FE})} has a polynomial solution without
multiple roots, then all solutions to this equation are polynomials.
\hfill $\triangleleft$
\end{prop}

\subsection{The proof of the left-hand inequality of  {\rm (\ref{s3})}}\label{s32}
Denote
\begin{equation}\label{Z} 
Z=\left\{\ \z=(z_1,\dots,z_n)\in\C^n\ \vert\  z_i\neq  z_j\,,\ 
1\leq i<j\leq n\ \right\}\,.
\end{equation}
Let $\z\in Z$ be fixed and $W(x)=W(x;\z,\m)$ be given by  (\ref{W}).
Define a symmetric function in $k$ complex
variables $\tb=(t_1,\dots, t_k)$,
\begin{equation}\label{Phi}
\Phi(\tb)=\Phi_{k,n}(\tb;\z,\m)\,=\,\frac{\prod_{j< l}(t_j-t_l)^2}
{\prod_{i=1}^k W(t_i)}\,,
\end{equation}
in the domain 
\begin{equation}\label{T}
T=T_{k,n}(\tb;\z)=\left\{\tb\in \C^k\ \vert \ t_i\neq t_l\,,\ t_i\neq  z_j\,,\
1\leq i,  l\leq k\,,\ 1\leq j\leq n\   \right\}.
\end{equation}
The function $\Phi(\tb)$ is called {\it the master function}.

\medskip
{\it A critical point} of  the function  $\Phi(\tb)$ is a 
point $\tb^0\in T$ such that 
$$
\frac{\partial\Phi}{\partial t_i}(\tb^0)=0\,,\quad i=1,\dots,k\,.
$$

\medskip
For any plane in $G_2(\Poly)$ with the Wronskian of degree $M$,  the sum
of the order and of the degree is $M+1$.
Under the assumptions (\ref{A}), we have $0 < M+1-d < d$.

\medskip
Heine and Stieltjes, in their studies of  second order linear differential 
equations with polynomial coefficients and a polynomial solution of 
a prescribed degree, arrived at the  result which can be formulated 
as follows. 

\begin{prop}\label{pW} {\rm (Cf. \cite{Sz},  Ch.~6.8)}\quad
Let $\z\in Z$ be fixed, $n, d, m_1,\dots, m_n$ be positive integers 
satisfying  {\rm (\ref{A})}, $M=m_1+\dots +m_n$, $k=M+1-d$.

If  $\tb^0$  is a critical point of  the master function {\rm (\ref{Phi})},
then $f(x)=(x-t^0_1)\dots(x-t^0_k)$ is  a  solution to  the equation  
{\rm (\ref{FE})}. 

Conversely,  if  $\tb^0\in T$ and if $f(x)=(x-t^0_1)\dots (x-t^0_k)$ is 
a  solution to the equation  {\rm (\ref{FE})}, then $t^0$ is a critical point 
of  the master function {\rm (\ref{Phi})}.
\hfill $\triangleleft$   
\end{prop}

The symmetric group  $S^k$ acts on $T$ permuting $t_1,\dots, t_k$,
and the action preserves the critical set of  the master function.
Say that  a generic plane $V\in G_2(\Poly)$  of order $k>0$ is   
{\it nondegenerate}, if the polynomials of   degree $k$ 
in $V$  do not have multiple roots. The nondegenerate planes
correspond to the classes of rational functions which can be 
represented by the ratio of polynomials without multiple roots.
Propositions~\ref{pp} and  \ref{pW} imply the following result. 

\begin{cor}\label{c3}
Let $\z\in Z$ be fixed, $W(x)=W(x;\z,\m)$ be given by  
{\rm (\ref{W}), (\ref{A})}, $k=M+1-d$, and $\Phi(\tb)=
\Phi_{k,n}(\tb;\z,\m)$ be given by {\rm (\ref{Phi})}.
There is a one-to-one correspondence between the orbits
of critical points of the master function $\Phi(\tb)$ 
and the  nondegenerate planes of degree $d$  with the Wronskian 
$W(x)$.
\hfill $\triangleleft$   
\end{cor}
Any  nondegenerate plane is generic, and  in order  to get
the left-hand inequality of (\ref{s3}) it remains to  use Proposition~\ref{p1}.
\hfill $\triangleleft$   

\medskip
\noindent
{\bf Remarks:}\quad\
\begin{itemize}
\item
For $k=1$, the master function is
$\Phi_{n,1}(t;\z,\m)=W^{-1}(t)$, and
the critical points  are exactly those of the critical points
of  $W(x)$ which have a non-zero critical value.
\item
According to  Proposition \ref{p2}, if $f(x)=(x-t_1^0)\dots(x-t_k^0)$ is
a polynomial without multiple roots in a nondegenerate plane 
of degree $d$, of order $k=M+1-d$ and with the Wronskian $W(x)$, 
then $\tb^0$ is a solution to the system
$$
\frac{W'(t_i)}{W(t_i)}\,=\,\frac{f''(t_i)}{f'(t_i)}\,,\quad  i=1,\dots,k.
$$ 
This is exactly the critical point system of the master function.
\item
After the first version of this paper was posted on arXiv, the author
received from A.~Eremenko (unpublished) notes where the
connection between classes of rational functions and Fuchsian equations
with trivial monodromy, in the case  $n=2d-2$ and $m_1=\dots=m_n=1$, 
has been established.
\end{itemize}

\section{A lower bound: asymptotic solutions to the $sl_2$ KZ equation}\label{S4}

\subsection{The $sl_2$ KZ equation {\rm (\cite{RV})}}\label{s41}
Consider the Lie algebra $sl_2=sl_2(\C)$ with the standard
 basis $e, f, h$ such that $[e,f]=h,\  [h,e]=2e,\ [h,f]=-2f\,.$
     
\medskip
For  positive integers  $m_1,\dots, m_n$, $n\geq 2$,  let  $L_{m_j}$  be 
the irreducible $sl_2$-module with  highest weight $m_j$, and
$L=L_{m_1}\otimes \dots \otimes  L_{m_n}$ the tensor product.

\medskip
For the Casimir element,                                                    
$$
\Omega=e\otimes f+f\otimes e+\frac12 h\otimes h \in sl_2\otimes sl_2,
$$
and for $1\leq i<j\leq n$,  let $\Omega_{ij}:L\to L$ be  
the operator which acts  as $\Omega$ on $i$-th and $j$-th factors of $L$
and as the identity on all others. For  $\z\in Z$, see (\ref{Z}),
define operators  $H_1(\z),\dots, H_n(\z)$ on $L$ as follows,
\begin{equation}
\label{H}
H_i(\z)=\sum_{j\neq i}\frac{\Omega_{ij}}{z_i-z_j}\,, \quad  i=1,\dots, n.
\end{equation}

\medskip
{\it The  Knizhnik--Zamolodchikov (KZ) equation} on a function 
$u: Z\to L$  is the system of partial differential equations
\begin{eqnarray}
\label{KZ}
\kappa \frac{\partial u}{\partial z_i}=H_i(\z)u(\z)\,,\quad   i=1,\dots, n,
\end{eqnarray}         
where $\kappa$ is a parameter. This equation appeared first 
in Wess--Zumino models of conformal field theory, \cite{KZ}.

\subsection{Subspaces of singular vectors}\label{s42}
Let $k$ be a nonnegative integer, $k\leq  M/2$ (as before,
$M=m_1+\dots+m_n$). Denote ${\Sing}_k L$ the subspace of 
singular vectors of weight  $M-2k$ in $L$,
\begin{equation}
\label{Sing}
{\Sing}_k L=\left\{\, w\in  L\ \vert\  ew=0\,,\ hw=(M-2k)w\, \right\}.
\end{equation}

\subsection{The hypergeometric solutions to the KZ equation}\label{s43}
The KZ equation preserves  ${\Sing}_kL $ for any $k$.
In \cite{ScV}, the hypergeometric solutions to the KZ equation
with values in ${\Sing}_k L$ were constructed.
Consider the function
 $$
\Psi(\tb,\z)=\prod_{1\leq i<j\leq n}(z_i-z_j)^{m_im_j/2}
\prod_{i=1}^k\prod_{l=1}^n (t_i-z_l)^{-m_l}
\prod_{1\leq i<j\leq k}(t_i-t_j)^2 \,
$$
defined on 
$$
\CC=\CC_{k,n}(\tb,\z)=\{\tb\in\C^k, \z\in\C^n\,\vert\,
t_i\neq t_l\,,\ z_j\neq z_p\,,\ t_i\neq z_j\,,\
1\leq i,l\leq k,\ 1\leq j,p\leq n\}.
$$
For the natural projection $\,\CC_{k,n}(\tb,\z)\to Z$,\ $\,(\tb,\z)\mapsto\z$,
the preimage of $\z$ is $T_{k,n}(\tb;\z)$, see (\ref{Z}), (\ref{T}).
In Sec.~7  of  \cite{ScV}, the solutions to the KZ equation   
(\ref{KZ})  with values in ${\Sing}_k L$ were constructed. They have the form
\begin{equation}\label{u}
u(\z)^{\gamma(\z)}=\int_{\gamma(\z)}\  \Psi^{1/\kappa}(\tb,\z)\, R(\tb,\z)\,
dt_1\wedge\dots\wedge dt_k\,,
\end{equation}
where   $\gamma(\z)$ is  a  $k$-cycle in $T=T_{k,n}(\tb;\z)$
continuously depending on $\z$, and  $R(\tb,\z)$ is  a certain rational 
function on $\CC$ with values in ${\Sing}_k L$.

\subsection{The proof of the right-hand inequality of {\rm (\ref{s3})}}\label{s44}
The series
$$
\phi(\z)=e^{I(\z)/\kappa}\left(\phi_0(\z)+\kappa \phi_1(\z)+\kappa^2 \phi_2(\z)+
\dots \right),
$$
where $I$ and $\phi_i$ are  functions on $Z$, provides
{\it an asymptotic solution to the equation}  (\ref{KZ}) if the substitution                             
of $\phi$ into (\ref{KZ}) gives $0$ in the
expansion into a formal power series in $\kappa$.   
If $\phi(\z)$ is  an asymptotic solution, then the  substitution  into   
(\ref{KZ}) shows that  $\phi_0(\z)\in {\Sing}_k L$  is  an eigenvector 
of  $H_j(\z)$  for any $1\leq j\leq n$.

\medskip
Asymptotic solutions can be produced by taking the limit as $\kappa \to 0$
of $u(\z)^{\gamma(\z)}$ given by (\ref{u}).
According to the steepest descend method,  
the leading terms $\phi_0(\z)$ are determined by critical points of the function
$\Psi^{1/\kappa}(\tb,\z)$ with respect to $\tb$. In order to study the 
critical points, one can clearly replace $\Psi^{1/\kappa}(\tb,\z)$ with 
the function  
$$
\Phi(\tb)=\prod_{i=1}^k\prod_{l=1}^n (t_i-z_l)^{-m_l}
\prod_{1\leq i<j\leq k}(t_i-t_j)^2
$$
defined on $T_{k,n}(\tb;\z)$ given by (\ref{T}). 
This is exactly the master function (\ref{Phi}).

\begin{prop}\label{pRV}{\rm (Theorem~9.9 of \cite{RV}, 
see also Theorem~8 of \cite{SV})}
\begin{itemize}
\item
Any orbit of  nondegenerate critical point $\tb^0(\z)$  of  the master function 
$\Phi(\tb)$
defines a common eigenvector  $w(\tb^0,\z)$ of  operators 
$H_1(\z),\dots, H_n(\z)$.
\item
For generic $\z\in Z$, the eigenvectors $w(\tb^0,\z)$ generate the space ${\Sing}_k L\,$. 
\hfill $\triangleleft$   
\end{itemize}
\end{prop}

Thus the number of orbits of critical points of the master function  
{\rm (\ref{Phi})} is at least $\dim{\Sing}_k L$. \hfill $\triangleleft$   

\section{Comments}\label{S5}
\subsection{A Hurwitz problem}\label{s51}
The question about the number of  rational maps
with prescribed singularities is a classical one.

\medskip
In 1891, A.~Hurwitz published a formula giving the number
of some special factorizations of a permutation into
transpositions, \cite{H}. This formula immediately
yields {\it the number of  topologically non-equivalent 
rational functions on $\C^1$ with  fixed orders 
of poles and fixed critical values, in the case when 
all critical points are nondegenerate and all critical 
values are different}. The topological 
equivalence of two functions means that one turns into 
the other by a linear-fractional transformation of the domain.

\medskip
Concerning the proof of the formula, there was a kind of 
discussion whether A.~Hurwitz
had solved it completely, or his proof contains a gap.
Since then, many papers on the subject were appeared.
Let us mention the paper of  I.P.~Goulden and D.M.~Jakson \cite{GJ} 
where  the formula was rediscovered in its full, cf. \cite{M}. 

\medskip
The question posed in the Introduction differs from the
Hurwitz's one. Our setting follows that of \cite{G}, 
\cite{EG}, cf. \cite{EH}. Firstly, not critical values  
but critical points are fixed. Secondly, a degree is
fixed as well. Thirdly, the equivalence relation is 
defined by linear-fractional transformations of the 
target, not the domain.

\subsection{Multiplicities of  $sl_2$-representations}\label{s52}
Formula (\ref{dim})  gives  the multiplicity of
$L_{M-2k}$ in the decomposition of the tensor product
$L_{m_1}\,\otimes\, \dots\, \otimes \, L_{m_n}$ into the direct sum
of irreducible  $sl_2$-modules.
The proof of the formula given in \cite{SV}  is  short. It
involves in  the Inclusion-Exclusion principle and  elementary
$sl_2$-representation theory.  However we were unable to find
this formula in the literature as on representation theory as
on the Schubert calculus.

\subsection {Degrees of Wronski maps and Catalan numbers}\label{s53}
The Wronski map sends the big cell of $G_2({\Poly}_d)$ 
(the planes of  degree $d$ and order  $d-1$) to the big cell of the projective
space (the monic polynomials of degree exactly $2d-2$).
Consider the restriction of the Wronski map to  the big cell.
The cardinality of a generic preimage is called 
{\it the degree of the Wronski map},  \cite{EG}.
 
\medskip
The main result of \cite{G} says that over $\C$ the degree of
the Wronski map equals the Catalan number $C_d=\frac 1d{2d-2\choose
d-1}$. In our notation, $C_d=\sharp(d,2d-2;\1;\z)$,  
where $\1=(1,\dots,1)$. The Catalan number gives also an upper bound 
for the number of the classes of rational maps of the type $(d,n;\m;\z)$ with
$m_1+\dots+m_n\leq 2d-2$ (\cite{G}).
In \cite{EG}, upper and lower estimates for the number
of classes of {\it real} rational maps of degree $d$ in terms of the
Catalan numbers were established.

\medskip
The Catalan numbers appear as the answer to many combinatorial
problems (see for example  \cite{St}, Ex.~6.19 and 6.25).
We get  that
{\it the $(k+1)$-th Catalan number is the multiplicity
of $L_0$ in $L_1^{\otimes 2k}$.}
 This interpretation leads to the following conclusion.

\begin{prop}
The generating function for the multiplicity $M_k$ of the trivial
$sl_2$ representation $L_0$ in $L_1^{\otimes k}$ is
$$
M(\tau)=\frac{1-2\tau^2-\sqrt{1-4\tau^2}}{2\tau^2}\,=
\,\sum_{k=1}^{\infty}M_k\tau^k\,,
\quad \tau\rightarrow 0.
$$
\end{prop}
 Indeed, the generating function for the Catalan numbers $C_k$
is as follows,  \cite{St},
$$
C(\tau)=\frac{1-\sqrt{1-4\tau}}2\,=\,\sum_{k=1}^{\infty}C_k \tau^k\,,
\quad \tau\rightarrow 0,
$$
and the multiplicity of $L_0$ in  $L_1^{\otimes 2k-1}$
is obviously zero. Therefore we have
$$
M(\tau)=\sum_{k=1}^{\infty}C_{k+1}\tau^{2k}=\frac{C(\tau^2)}{\tau^2}-1\,.
$$
This gives the statement.    \hfill $\triangleleft$   

\subsection{The Schubert calculus and Fuchsian equations 
with only polynomial solutions}\label{s54}
One of the main results of \cite{SV} is as follows.
\begin{prop} {\rm (Theorem 3 of \cite{SV})} \quad
Let $m_1, \dots, m_n,$ and $k$ be nonnegative integers, 
$n\geq 2$, $M=m_1+\dots+m_n$ and $k\leq M/2$. 
Then for generic $\z\in Z$ there exist exactly 
$$
\sum_{q=0}^{n} (-1)^q
\sum_{1\leq i_1<\dots <i_q\leq n}
{k+n-2-m_{i_1}-\dots -m_{i_q}-q\choose n-2}
$$                        
polynomials $H(x)$ of degree not greater than $n-2$ such that all
solutions to the equation {\rm (\ref{FE})} are polynomials and the degree
of the generic solution equals $M+1-k$.   \hfill $\triangleleft$    
\end{prop}

The proof  given  in \cite{SV} is rather complicated. Using of  
the Schubert calculus leads to another, ``less technical'', proof,
in view of Propositions~\ref{p1}, \ref{pF}.

\medskip
The theory of Fuchsian differential equations of order $p$ with only 
polynomial solutions  is closely related to the Schubert calculus 
in the Grassmannian of $p$-dimensional subspaces of the vector space
of complex polynomials. If the generic solution of such an equation  
has degree $d$, then the solution space $V$ belongs to $G_p({\Poly}_d)$.
Let polynomials  $P_1(x),\dots, P_p(x)$ form a basis of
$V\in G_p({\Poly}_d)$. Then the equation $E_V(u)=0$ with respect to unknown
function $u$, where
$$
E_V(u)=
\det \left(\begin{array}{cccc} u(x)\ \ &P_1(x)\   \ &\dots\   \ &P_p(x)\ \ \\
                          u'(x)\ \ &P_1'(x)\ \  &\dots\ \ &P_p'(x)\    \\
                           \dots\ \  &\dots\ \   &\dots\ \  &\dots\ \\
                          u^{(p)}(x)\   &P_1^{(p}(x)\ \ &\dots\ \ &P_p^{(p)}(x)\\
                           \end{array}\right),
$$
is Fuchsian and its solution space is $V$.
{\it The Wronskian $W_V(x)$ of $V$}  is defined as a monic polynomial
which is proportional to the Wronskian  of some (and hence, any) basis of  $V$,
cf. \cite{EG}. The Wronskian of the equation $E_V(u)=0$ is proportional 
to  $W_V(x)$, and the regular singular points of the equation are 
the roots of  $W_V(x)$ (``the root at infinity" has
multiplicity $p(d+1-p)-\deg W_V$).
For any $z\in\C\cup\infty$, the set of exponents at $z$ defines 
a Schubert variety $\Omega_{w(z)}(\F_\bullet(z))$  containing $V$  
(CH.~6 of \cite{GH}, cf.~ Sec.~\ref{s25}).
In particular, $z$ is an ordinary point of the equation if and only if
$w(z)=(0,\dots,0)$. The intersection of the Schubert varieties 
corresponding to all regular singular points consists of a finite number
of elements. Thus the number of  Fuchsian differential
equations of order $p$ with only polynomial solutions, having
given regular singular points,  given exponents at these points and the 
degree $d$  generic solution, is  bounded from above by the intersection
number of the corresponding Schubert classes.

\subsection{The Schubert calculus and the $sl_p$ KZ equations}\label{s55}
In our case, the theory of  KZ equations helps to solve a problem
from enumerative algebraic geometry.  On the other hand,
the Schubert calculus seems to be a useful tool in the theory of 
KZ equations. In particular, the function which appears in the
hypergeometric solutions to the $sl_p$ KZ equation   
plays the role of the master function for the Wronski map from the 
Grassmannian of $p$-dimensional  planes, \cite{Ga}, \cite{S}. 
The Schubert calculus provides an upper bound for the number of orbits 
of critical points as the intersection number of appropriate Schubert classes. 
For the tensor product of  symmetric powers of the standard $sl_p$-module,
a lower bound can be obtained,  combining   the theory of Fuchsian equations 
and  results  on KZ equations, as the dimension of a relevant
subspace of singular vectors, \cite{S}. 
The coincidence of these upper and lower bounds implies  non-trivial 
results  both  in the  theory of  KZ equations  (the proof of the conjecture 
on Bethe vectors in Gaudin model) and in the enumerative algebraic geometry 
(the proof of the transversality claim for Schubert varieties).

\newpage
\

\end{document}